\documentclass[12pt]{amsart}
\usepackage{amssymb,amsmath,amsfonts,latexsym}
\usepackage{bm}
\newcommand{\newsection}[1]{\setcounter{equation}{0} \section{#1}}
\newcommand{\bea}{\begin{eqnarray}}
\newcommand{\eea}{\end{eqnarray}}

\newcommand{\clb}{\mathcal{B}}

\newcommand{\cld}{\mathcal{D}}
\newcommand{\cle}{\mathcal{E}}

\newcommand{\clh}{\mathcal{H}}
\newcommand{\clk}{\mathcal{K}}
\newcommand{\cll}{\mathcal{L}}
\newcommand{\clm}{\mathcal{M}}

\newcommand{\clo}{\mathcal{O}}

\newcommand{\clq}{\mathcal{Q}}
\newcommand{\clr}{\mathcal{R}}
\newcommand{\cls}{\mathcal{S}}

\newcommand\bv{\bm{v}}
\newcommand\T{\bm{T}}
\newcommand{\D}{\mathbb{D}}
\newcommand{\N}{\mathbb{N}}
\newcommand{\C}{\mathbb{C}}
\newcommand{\z}{\bm{z}}
\newcommand{\w}{\bm{w}}

\newcommand{\raro}{\rightarrow}
\newcommand{\Raro}{\Rightarrow}

\def \qed {\hfill \vrule height6pt width 6pt depth 0pt}
\def\textmatrix#1&#2\\#3&#4\\{\bigl({#1 \atop #3}\ {#2 \atop #4}\bigr)}
\def\dispmatrix#1&#2\\#3&#4\\{\left({#1 \atop #3}\ {#2 \atop #4}\right)}
\newcommand{\be}{\begin{equation}}
\newcommand{\ee}{\end{equation}}
\newcommand{\ben}{\begin{eqnarray*}}
\newcommand{\een}{\end{eqnarray*}}

\newcommand{\NI}{\noindent}

\newcommand{\bi}{\begin{itemize}}
\newcommand{\ei}{\end{itemize}}

\newcommand\m{\bm{m}}

\newcommand{\DS}{\mathcal{D}_{\bm{m},\T}}

\newtheorem{Theorem}{\sc Theorem}[section]
\newtheorem{Lemma}[Theorem]{\sc Lemma}
\newtheorem{Proposition}[Theorem]{\sc Proposition}
\newtheorem{Corollary}[Theorem]{\sc Corollary}
\newtheorem{Definition}[Theorem]{\sc Definition}
\newtheorem{Example}[Theorem]{\sc Example}
\newtheorem{Remark}[Theorem]{\sc Remark}
\newtheorem{Note}[Theorem]{\sc Note}
\newtheorem{Question}{\sc Question}
\newtheorem{ass}[Theorem]{\sc Assumption}
\newcommand{\bt}{\begin{Theorem}}
\def\beginlem{\begin{Lemma}}
\def\beginprop{\begin{Proposition}}
\def\begincor{\begin{Corollary}}
\def\begindef{\begin{Definition}}
\def\beginexamp{\begin{Example}}
\def\beginrem{\begin{Remark}}
\def\beginq{\begin{Question}}
\def\beginass{\begin{ass}}
\def\beginnote{\begin{Note}}
\newcommand{\et}{\end{Theorem}}
\def\endlem{\end{Lemma}}
\def\endprop{\end{Proposition}}
\def\endcor{\end{Corollary}}
\def\enddef{\end{Definition}}
\def\endexamp{\end{Example}}
\def\endrem{\end{Remark}}
\def\endq{\end{Question}}
\def\endass{\end{ass}}
\def\endnote{\end{Note}}

\begin{document}

\title[Operator Positivity and Analytic Models]{Operator Positivity and Analytic Models of Commuting Tuples of Operators}

\author[Monojit Bhattacharjee]{Monojit Bhattacharjee}
\address{Indian Institute of Science, Department of Mathematics, Bangalore, 560012, India}
\email{monojit12@math.iisc.ernet.in}

\author[Jaydeb Sarkar]{Jaydeb Sarkar}
\address{Indian Statistical Institute, Statistics and Mathematics Unit, 8th Mile, Mysore Road, Bangalore, 560059, India}
\email{jay@isibang.ac.in, jaydeb@gmail.com}

\subjclass[2000]{47A13, 47A15, 47A20, 47A45, 47A80, 47B32, 47B38}


\keywords{Weighted Bergman spaces, hypercontractions, multipliers,
reproducing kernel Hilbert spaces, invariant subspaces}

\begin{abstract}
We study analytic models of operators of class $C_{\cdot 0}$ with
natural positivity assumptions. In particular, we prove that for an
$m$-hypercontraction $T \in C_{\cdot 0}$ on a Hilbert space $\clh$,
there exists two Hilbert spaces $\cle$ and $\cle_*$ and a partially
isometric multiplier $\theta \in \clm(H^2(\cle), A^2_m(\cle_*))$
such that
\[
\clh \cong \clq_{\theta} = A^2_m(\cle_*) \ominus \theta H^2(\cle), \quad \quad \mbox{and} \quad \quad T \cong
P_{\clq_{\theta}} M_z|_{\clq_{\theta}},
\]
where $A^2_m(\cle_*)$ is the $\cle_*$-valued weighted Bergman space
and $H^2(\cle)$ is the $\cle$-valued Hardy space over the unit disc
$\D$. We then proceed to study and develop analytic models for
doubly commuting $n$-tuples of operators and investigate their
applications to joint shift co-invariant subspaces of reproducing
kernel Hilbert spaces over polydisc. In particular, we completely
analyze doubly commuting quotient modules of a large class of
reproducing kernel Hilbert modules, in the sense of Arazy and
Englis, over the unit polydisc $\D^n$.
\end{abstract}

\maketitle

\section*{Notation}

\begin{list}{\quad}{}
\item $\mathbb{N}$ \quad \quad \; Set of all natural numbers
including 0.

\item $n$ \quad \quad \; Natural number $n \geq 2$.
\item $\mathbb{N}^n$ \quad \quad $\{\bm{k} = (k_1, \ldots, k_n) : k_i \in \mathbb{N}, i = 1,
\ldots, n\}$.
\item $\bm{z}$ \quad \quad \; $(z_1, \ldots, z_n) \in \mathbb{C}^n$.
\item $\bm{z}^{\bm{k}}$ \quad \quad \,$z_1^{k_1}\cdots
z_n^{k_n}$.
\item $\T$ \quad \quad \; $n$-tuple of commuting operators $(T_1, \ldots, T_n)$.
\item $\T^{\bm{k}}$ \quad \quad $T_1^{k_1} \cdots
T_n^{k_n}$.
\item $\mathbb{D}^n$  \quad \quad Open unit polydisc $\{\z : |z_i|
<1\}$.
\end{list}

For a closed subspace $\cls$ of a Hilbert space
$\clh$, we denote by $P_{\cls}$ the orthogonal projection of $\clh$
onto $\cls$. We shall denote the space of all bounded linear operators on $\clh$ by $\clb(\clh)$.

\newsection{Introduction}

The Sz.-Nagy and Foias analytic model theory for contractions on
Hilbert spaces is a powerful tool for studying operators on Hilbert
spaces and holomorphic function spaces on the open unit disc $\D$ in
$\C$. It says that if $T$ is a contraction (that is, $I - T T^* \geq
0$) on a Hilbert space and in $C_{\cdot 0}$ class (that is, $T^{*l}
\raro 0$ as $l \raro \infty$ in strong operator topology) then $T^*$
is unitarily equivalent to the restriction of the backward shift
$M_z^*$ on a vector-valued Hardy space to a $M_z^*$-invariant
subspace. More precisely, there exists a coefficient Hilbert space
$\cle_*$ and a $M_z^*$-invariant closed subspace $\clq$ of
$\cle_*$-valued Hardy space $H^2(\cle_*)$ such that \[T \cong
P_{\clq} M_z|_{\clq}.\] Moreover, there exists a Hilbert space
$\cle$ and a $\clb(\cle, \cle_*)$-valued inner multiplier $\theta_T
\in H^\infty_{\clb(\cle, \cle_*)}(\D)$, also known as the
characteristic function of $T$ (see \cite{NF}), such that
\[
\clq = H^2(\cle_*)/\theta_T H^2(\cle).
\]

In \cite{Ag}, Agler introduced and studied the theory
hypercontraction operators from operator positivity point of view.
He showed that the vector-valued Hardy space in the dilation space
of a contraction can be replaced by a vector-valued weighted Bergman
space if the contractivity assumption on the operator is replaced by
a weighted Bergman type positivity. Later, Muller and Vasilescu
\cite{MV}, Curto and Vasilescu \cite{CV}, Ambrozie and Timotin
\cite{AT1, AT2}, Arazy, Englis and Muller
\cite{AEM} and Arazy and Englis \cite{AE} extended these ideas to a
more general class of operators. This viewpoint has proved to be
extremely fruitful in studying commuting tuples of operators.

The purpose of this paper is to explore how one might do analytic
model theory for a general class of operators and commuting tuples
of operators. In particular, we associate a partially isometric
multiplier with every operator satisfying weighted Bergman-type
positivity condition (see Theorem \ref{model-1}). Another basic
result in this direction is the following analytic model: Let $\T =
(T_1, \ldots, T_n)$ be a doubly commuting tuple of pure operators on
a Hilbert space $\clh$ (that is, $T_i \in C_{\cdot 0}$, $T_i T_j =
T_j T_i$ and $T_p T_q^* = T_q^* T_p$ for all $i, j = 1, \ldots, n$,
and $1 \leq p < q \leq n$). Then $(T_1^*, \ldots, T_n^*)$ is joint
unitarily equivalent to the restriction of $(M_{z_1}^*, \ldots,
M_{z_n}^*)$ to a joint invariant subspace of a vector-valued
weighted Bergman space if and only if $T$ satisfies (joint) weighted
Bergman type positivity. Moreover, in this case, the orthocomplement
of the co-invariant subspace of the weighted Bergman space is of
``Beurling-lax-Halmos'' type (see Theorems \ref{MT}).

Although our method works for more general cases (see Section 7),
for simplicity we will restrict our discussion to hypercontractions
(see Section 2).

Here is a brief description of the paper. In Section 2 we set up
notations, recall some basic notions from the theory of
hypercontractions and construct an analytic structure on the model
space. Our main tool here is the Agler dilation theorem for
hypercontractions \cite{Ag} combined with a Beurling-Lax-Halmos type
representations of shift invariant subspaces of analytic reproducing
kernel Hilbert spaces (\cite{BV1}, \cite{JS}). In Section 3 we
discuss a dilation theory for a class of doubly commuting operator
tuples satisfying weighted Bergman type positivity condition.  In
Section 4 we formulate a version of Sz.-Nagy and Foias analytic
model for doubly commuting tuple of hypercontractions. In Section 5,
we analyze doubly commuting quotient modules of scalar valued
weighted Bergman space. Finally in section 6 we study
$K$-contractive tuples of operators in the spirit of Arazy and
Englis \cite{AE}.

\section{Functional Models for Hypercontractions}

The main purpose of this section is to study and to develop an
analytic functional model for the class of hypercontractions on
Hilbert spaces.

We first recall the definition of weighted Bergman spaces and review
the construction of the dilation maps for hypercontractions. We
refer the reader to Agler's paper \cite{Ag} for more
details.

The weighted Bergman kernel on the open unit disc $\D$ with weight
$\alpha >0$ is, by definition, the kernel function:
\[
B_{\alpha}(z, w) = (1 - z \bar{w})^{-\alpha} \quad \quad (z, w
\in \D).
\]
For each $\alpha >0$, we let $A^2_{\alpha}$, denote the weighted
Bergman space corresponding to the kernel $B_{\alpha}$. For any
Hilbert space $\cle$, the $\cle$-valued weighted Bergman space
$A^2_{\alpha}(\cle)$ with reproducing kernel $(z, w) \in \D \times
\D \mapsto B_{\alpha}(z, w) I_{\cle}$ can canonically be identified
with the Hilbert space tensor product $A^2_{\alpha} \otimes \cle$.
In order to simplify notation, we often identify $A^2_{\alpha}
\otimes \cle$ with $A^2_{\alpha}(\cle)$. It also follows that
$\{B_{\alpha}(\cdot, w) \eta: w \in \D, \eta \in \cle\}$ is a total
set in $A^2_{\alpha}(\cle)$ and $\langle f, B_{\alpha}(\cdot, w)
\eta \rangle_{A^2_{\alpha}(\cle)} = \langle f(w), \eta
\rangle_{\cle}$ where $f \in A^2_{\alpha}(\cle)$, $w \in \D$ and
$\eta \in \cle$. Moreover, it is easy to see that the shift operator
$M_z$ on $A^2_{\alpha}(\cle)$, $\alpha \geq 1$, is a $C_{\cdot
0}$-contraction where
\[
(M_z f)(w) = w f(w) \quad \quad \quad (f \in A^2_{\alpha}(\cle), w
\in \D).
\]

In the following discussion, we shall mostly use weighted Bergman
spaces with only integer weights. Let us point out an important
special case: $A^2_1 = H^2$, the Hardy space over $\D$.

For a multi-index $\m = (m_1, \ldots, m_n) \in \N^n$ we denote the
corresponding weighted Bergman space on $\D^n$ by $A^2_{\m}$. The
weighted Bergman kernel on $\D^n$ with weight $\m$ is, by
definition, the reproducing kernel function
\[
B_{\m}(\z, \w) = \prod_{i=1}^n B_{m_i}(z_i, w_i) = \prod_{i=1}^n (1
- z_i \bar{w_i})^{-m_i} \quad \quad (\z, \w \in \D^n).
\]
For each $\w \in \D^n$, we denote by $B_{\m}(\cdot, \w)$ the kernel
function at $\w$, where
\[
(B_{\m}(\cdot, \w))(\z) = B_{\m}(\z, \w) \quad \quad \quad (\z \in
\D^n).
\]

\NI\textsf{Convention:} Let $p(\z, \w) = \sum_{\bm{p}, \bm{q} \in
\N^n} a_{\bm{p} \bm{q}} \z^{\bm{p}} \bar{\w}^{\bm{q}}$ be a
polynomial in $\{z_1, \ldots, z_n\}$ and $\{\bar{w}_1, \ldots,
\bar{w}_n\}$. For a commuting tuple of bounded linear operators $\T =
(T_1, \ldots, T_n)$ on a Hilbert space $\clh$ (that is, $T_i T_j =
T_j T_i$ for all $i, j = 1, \ldots, n$) we denote by $p(\z, \w)(\T,
\T^*)$ the corresponding hereditary functional calculus in the sense
of Agler \cite{Ag}, that is,
\begin{equation}\label{hered}p(z, w)(\T, \T^*) = \sum_{\bm{p}, \bm{q} \in \N^n} a_{\bm{p}
\bm{q}} \T^{\bm{p}} \T^{* \bm{q}},
\end{equation}
where $\T^{\bm{k}} = T^{k_1} \cdots T^{k_n}_n$ and $\T^{*\bm{k}} =
T_1^{*k-1} \cdots T_n^{*k_n}$ for all $\bm{k} =(k_1, \ldots, k_n)
\in \mathbb{N}^n$.

\begin{Definition}
A bounded linear operator $T$ on $\clh$ is said to be
$B_m$-contractive (or $T$ is a $B_m$-contraction) if $T$ is in
$C_{\cdot 0}$ class and \[\begin{split}B_m^{-1}(z, w)(T, T^*) & =
\big(\sum_{k=0}^m (-1)^k {m\choose k} {z}^k \bar{w}^{*k}\big)(T, T^*) \\
&= \sum_{k=0}^m (-1)^k {m\choose k} T^k T^{*k} \geq 0.\end{split}\]
\end{Definition}

We also recall that $T \in \clb(\clh)$ is a
\textit{hypercontraction} of order $m$ \cite{Ag} if
\[B_p^{-1}(z, w)(T, T^*) \geq 0,\]holds for all $1 \leq p \leq m$.

Now let $T$ be a $B_m$-contraction on $\clh$. Since $T \in C_{\cdot
0}$, it follows from Lemma 2.11 in \cite{Ag} that
\[
B_p^{-1}(z, w)(T, T^*) \geq 0 \quad \quad (1 \leq p \leq m),
\]
that is, $T$ is a hypercontraction of order $m$. In other wards,
these two notions coincide for $C_{\cdot 0}$ class of operators and
hence, we will restrict our considerations for $B_m$-contractions.

The \textit{defect operator} and the \textit{defect space} of a
$B_m$-contraction $T \in \clb(\clh)$ are defined by
\begin{equation}\label{defect}
D_{m, T} = {\Big(B^{-1}_m(z, w)(T, T^*)\Big)}^{\frac{1}{2}}, \quad
\cld_{m,T} = \overline{\mbox{ran}}  D_{m, T},
\end{equation}
respectively. Let us set
\begin{equation}\label{BzT}
B_m(z, T) = (I_{\clh} - z T^*)^{-m}\quad  \quad \quad (z \in \D),
\end{equation}
and
\[
(\bm{v}_{m, T} f)(z) = D_{m,T} B_m(z, T) f = D_{m,T} (I_{\clh} - z
T^*)^{-m} f \quad \quad \quad (f \in \clh, z \in \D).
\]
Then $\bm{v}_{m, T} : \clh
\raro A^2_m (\cld_{m, T})$ is a bounded linear operator,
\[
\bm{v}_{m, T} T^* = M_z^*\bm{v}_{m, T},
\]
and
\[
\bv_{m, T}^*\Big(B_m(\cdot, w) \eta \Big) =
B_m(w, T)^* D_{m,T} \eta = (I_{\clh} - \bar{w} T)^{-m} D_{m,T}
\eta,
\]
for all $w \in \D$ and $\eta \in \cld_{m,T}$. This and the
definition of $\bv_{m,T}$ together imply
\begin{equation}\label{vv*}
\Big(\bv_{m, T} \bv_{m, T}^* (B_m(\cdot, w) \eta )\Big)(z) = D_{m,
T} B_m(z, T) B_m(w, T)^* D_{m,T}\eta \quad  \quad \quad (z \in \D),
\end{equation}
for all $w \in \D$ and $\eta \in \cld_{m, T}$. Furthermore, if $T$
is a $B_m$-contraction then $\bv_{m, T}$ is an isometry and hence a
dilation of $T$ (see Agler \cite{Ag}).

\begin{Theorem}\label{Agler}\textsf{(Agler )}
Let $T \in \clb(\clh)$ be a $B_m$-contraction. Then $T \cong
P_{\clq} M_z|_{\clq}$, for some $M_z^*$-invariant closed subspace
$\clq$ of $A^2_m(\cld_{m, T})$.
\end{Theorem}

We shall now introduce the notion of multipliers on weighted Bergman
spaces. Let $m_1, m_2$ be two natural numbers and $\cle_1, \cle_2$
be two Hilbert spaces. An operator valued holomorphic map $\Theta :
\mathbb{D} \raro \clb(\cle_1, \cle_2)$ is said to be a
\textit{multiplier} from $A^2_{m_1}(\cle_1)$ to $A^2_{m_2}(\cle_2)$
if $\Theta f \in A^2_{m_2}(\cle_2)$ for all $f \in
A^2_{m_1}(\cle_1)$. We denote the set of all multipliers from
$A^2_{m_1}(\cle_1)$ to $A^2_{m_2}(\cle_2)$ by
$\clm(A^2_{m_1}(\cle_1), A^2_{m_2}(\cle_2))$. We also use the
notation $M_{\Theta}$, for each $\Theta \in \clm(A^2_{m_1}(\cle_1),
A^2_{m_2}(\cle_2))$, to denote the multiplication operator
\[
M_{\Theta} f = \Theta f \quad \quad \quad (f \in
A^2_{m_1}(\cle_1)).
\]

\NI A multiplier $\Theta \in \clm(A^2_{m_1}(\cle_1),
A^2_{m_2}(\cle_2))$ is said to be a partially isometric multiplier
if $M_{\Theta}$ is a partially isometric operator from
$A^2_{m_1}(\cle_1)$ to $A^2_{m_2}(\cle_2)$.

Before proceeding, let us for the sake of completeness recall a
Beurling-Lax-Halmos type theorem for weighted Bergman shifts (see
\cite{BV1} and Theorem 2.3 in \cite{JS}) upon which much of our
discussion in this paper will rest.

\begin{Theorem}\label{MT-blh}
Let $\cls$ be a non-trivial closed subspace of $A^2_m(\cle_*)$. Then
$\cls$ is $M_z$-invariant if and only if there exists a Hilbert
space $\cle$ and a partially isometric multiplier $\theta \in
\clm(A^2_1(\cle), A^2_m(\cle_*))$ such that $\cls = \theta
A^2_1(\cle)$.
\end{Theorem}

We are now ready to present a functional model for the class of
$B_m$-contractions.

\begin{Theorem}\label{model-1}
Let $T \in \clb(\clh)$ be a $B_m$-contraction. Then there exists a
Hilbert space $\cle$ and a partially isometric multiplier $\theta
\in \clm(A^2_1(\cle), A^2_m(\cld_{m, T}))$ such that \[T \cong
P_{\clq_{\theta}} M_z|_{\clq_{\theta}},\]where $\clq_{\theta} =
A^2_m(\cld_{m, T}) \ominus \theta A^2_1(\cle)$.
\end{Theorem}
\NI\textsf{Proof.} At first, by virtue of Theorem \ref{Agler}, we
realize $T$ as $T \cong P_{\clq} M_z|_{\clq}$. Therefore, it only
remains to prove the existence of a partially isometric multiplier
$\theta$ such that $\clq = A^2_m(\cld_{m, T}) \ominus \theta
A^2_1(\cle)$.

\NI Note that since $\clq = \mbox{ran~}\bv_{m,T}$ is
$M_z^*$-invariant, $(\mbox{ran~}\bv_{m,T})^{\perp}$ is a
$M_z$-invariant closed subspace of $A^2_m(\cld_{m, T})$. Then,
applying Theorem \ref{MT-blh} to $(\mbox{ran~}\bv_{m,T})^{\perp}$,
we obtain a coefficient Hilbert space $\cle$ and a partially
isometric multiplier $\theta \in \clm(A^2_1(\cle), A^2_m(\cld_{m,
T}))$ such that
\[
(\mbox{ran~}\bv_{m,T})^{\perp} = \theta A^2_1(\cle),
\]that is,
\[
\clq = A^2_m(\cld_{m, T}) \ominus \theta A^2_1(\cle).
\]
This completes the proof. \qed

The following observation was pointed out to us by R. G. Douglas:
Let $T \in \clb(\clh)$ be a $B_m$-contraction. Then by Theorem
\ref{model-1}, we have
\[
H^2(\cle) \stackrel{M_{\theta}} \longrightarrow A^2_m(\cld_{m, T})
\stackrel{\pi} \longrightarrow \clh \longrightarrow 0,
\]
where $\pi = \bm{v}_{m,T}^*$ (recall that $A^2_1(\cle) =
H^2(\cle)$). Note that since $M_{\theta} M_z = M_z M_{\theta}$,
$\cls := \ker M_{\theta}$ is a $M_z$-invariant subspace of
$H^2(\cle)$. Then by Beurling-Lax-Halmos theorem there exists a
Hilbert space $\cle_*$ and an inner (or isometric) multiplier $\psi
\in H^\infty_{\clb(\cle_*, \cle)}(\D)$ such that $\cls = \psi
H^2(\cle_*)$. Consequently, we have a natural chain complex of
Hilbert spaces:\[0 \longrightarrow H^2(\cle_*) \stackrel{M_{\psi}}
\longrightarrow H^2(\cle) \stackrel{M_{\theta}} \longrightarrow
A^2_m(\cld_{m, T}) \stackrel{\pi} \longrightarrow \clh
\longrightarrow 0.\]

\newsection{Dilations of commuting hypercontractions}

In this section, we give a proof of the fact that a doubly commuting
tuple of hypercontractions can be dilated to the tuple of
multiplication operators on a suitable weighted Bergman space over
$\D^n$ (see \cite{AEM, AT1, AT2, AE}). We begin with a definition.

\begin{Definition}\label{Def1}
A commuting tuple of operators $\T = (T_1, \ldots, T_n)$ on $\clh$
is said to be $B_{\m}$-contractive if $T_i$ is a
$B_{m_i}$-contraction, $i = 1, \ldots, n$, and
\[B_{\m}^{-1}(\z, \w) (\T, \T^*) = \Big( \prod_{i=1}^n
B_{m_i}^{-1}(z_i, w_i) \Big)(\T, \T^*) \geq 0.\] We denote the
(joint-)defect operator and defect space of a $B_{\m}$-contraction
$\T$ as:
\[
D_{\m,\T} := (B_{\m}^{-1}(\z, \w) (\T, \T^*))^\frac{1}{2}, \quad \quad  \mbox{and} \quad\quad \cld_{\m,\T} =
\overline{\mbox{ran}} D_{\m,\T},
\]
respectively.
\end{Definition}

\textsf{For the rest of the paper we shall be dealing with a fixed
natural number $n \geq 2$, and a multi-index $\m = (m_1, \ldots,
m_n) \in \mathbb{N}^n$, $m_j \geq 1$, $j = 1, \ldots, n$.}

Let $\T$ be a doubly commuting $B_{\m}$-contractive tuple on $\clh$.
Then
\[
T_i (B^{-1}_{m_j}(z, w)(T_j, T_j^*)) = (B^{-1}_{m_j}(z, w)(T_j,
T_j^*)) T_i,
\]
for all $i \neq j$. This also implies that
\[
\Big(B_{m_i}^{-1}(z_i, w_i) (T_i, T_i^*) \Big)
\Big(B_{m_j}^{-1}(z_j, w_j) (T_j, T_j^*) \Big) =
\Big(B_{m_j}^{-1}(z_j, w_j) (T_j, T_j^*) \Big)
\Big(B_{m_i}^{-1}(z_i, w_i) (T_i, T_i^*) \Big).
\]

The above observations yield the following:

\begin{Lemma}\label{lemma1}
Let $\T$ be an $n$-tuple of doubly commuting operators on $\clh$ and
each $T_i$ is a $B_{m_i}$-contraction, $i = 1 , \ldots, n$. Then
$\cld_{m_j, T_j}$ is a $T_i$-reducing subspace of $\clh$ and
\[
T_i D_{m_j, T_j} = D_{m_j, T_j} T_i, \quad D_{m_j, T_j} D_{m_i, T_i}
= D_{m_i, T_i} D_{m_j, T_j},
\]
for all $i \neq j$. Moreover, $\T$ is a $B_{\m}$-contraction and
\[
D_{\m,\T} = \prod_{i=1}^n D_{m_i, T_i}.
\]
\end{Lemma}

Now, we shall construct, using induction, a dilation map for a
doubly commuting $B_{\m}$-contractive tuple $\T$ on $\clh$.

\textit{In what follows, for each $j \in \{2, \ldots, n\}$,
$\bm{m_j}$ denotes the $j$-tuple $(m_1, \ldots, m_j) \in
\mathbb{N}^j$ and $\bm{T_j}$ the doubly commuting
$B_{\bm{m_j}}$-contraction $(T_1, \ldots, T_j)$ on $\clh$.} For each
$j \in \{2, \ldots, n\}$, one checks easily that
$\cld_{\bm{m_{j-1}}, \bm{T_{j-1}}}$ is a $T_j$-reducing closed
subspace of $\clh$ and
\[
D_{T_j|_{\cld_{\bm{m_{j-1}}, \bm{T_{j-1}}}}} =
D_{T_j}|_{\cld_{\bm{m_{j-1}}, \bm{T_{j-1}}}}, \quad \quad \mbox{and}
\quad \quad \cld_{T_j|_{\cld_{\bm{m_{j-1}}, \bm{T_{j-1}}}}} =
\cld_{\bm{m_j}, \bm{T_j}},
\]
where $\cld_{\bm{m_1}, \bm{T_1}} = \cld_{m_1, T_1}$. Now we set $V_1
:= \bv_{m_1, T_1}: \clh \raro A^2_{m_1}(\cld_{m_1, T_1})$ and define
\[
V_2 := I_{A^2_{m_1}} \otimes \bv_{m_2, T_2|_{\cld_{m_1, T_1}}} :
A^2_{m_1} \otimes \cld_{m_1, T_1} \cong A^2_{m_1}(\cld_{m_1, T_1})
\raro A^2_{m_1} \otimes A^2_{m_2} \otimes \cld_{\bm{m_2}, \bm{T_2}}
\cong A^2_{\bm{m_2}} (\cld_{\bm{m_2}, \bm{T_2}}),
\]
where $\bv_{m_2, T_2|_{\cld_{m_1, T_1}}}: \cld_{m_1, T_1} \raro
A^2_{m_2} \otimes \cld_{\bm{m_2}, \bm{T_2}}$ is the dilation of
$T_2|_{\cld_{m_1, T_1}} \in \clb(\cld_{m_1, T_1})$. It follows that,
for all $l \in \N, h \in \cld_{m_1, T_1}$,
\[
\Big(V_2(z^l h)\Big)(z_1, z_2) = z_1^l (\bv_{m_2, T_2} h)(z_2).
\]
Continuing this way, one can construct $n$ bounded linear operators
$\{V_j\}_{j=2}^n$ defined by
\[
V_j = I_{A^2_{\bm{m_{j-1}}}} \otimes \bm{v}_{m_j,
T_j|_{\cld_{\bm{T_{j-1}}}}} : A^2_{\bm{m_{j-1}}}
(\cld_{\bm{m_{j-1}}, \bm{T_{j-1}}}) \longrightarrow A^2_{\bm{m_j}}
(\cld_{\bm{m_j}, \bm{T_j}}),
\]
where
\begin{equation}\label{V_j}
\Big(V_j(z_1^{k_1} \cdots z_{j-1}^{k_{j-1}} h)\Big)(z_1, \ldots,
z_j) = z_1^{k_1} \cdots z_{j-1}^{k_{j-1}} (\bm{v}_{m_j,
T_j|_{\cld_{\bm{m_{j-1}}, \bm{T_{j-1}}}}} h)(z_j),
\end{equation}
for all $h \in \cld_{\bm{m_{j-1}}, \bm{T_{j-1}}}$ and $j = 2,
\ldots, n$. Consequently, we have the following sequence of maps:
\[
0\longrightarrow \clh \stackrel{V_1} \longrightarrow
A^2_{m_1}(\cld_{m_1, T_1}) \stackrel{V_2} \longrightarrow
A^2_{\bm{m_2}}(\cld_{\bm{m_2}, \bm{T_2}})
\stackrel{V_3}\longrightarrow \cdots \stackrel{V_n}\longrightarrow
A^2_{\m}(\cld_{\bm{m}, \T}).
\]
Let us denote by $V_{\T}$ the compositions of $\{V_j\}_{j=1}^n$:
\begin{equation}\label{defn-VT}
V_{\T} := V_n \circ \cdots \circ V_2 \circ V_1 : \clh \raro
A^2_{\m}(\cld_{\m, \T}).
\end{equation}
Then $V_{\T} \in \clb(\clh, A^2_{\m}(\cld_{\m, \T}))$ is an
isometric dilation of $\T$:

\begin{Theorem}\label{dil1}
Let $\T$ be a doubly commuting $B_{\m}$-contractive tuple on $\clh$.
Then $V_{\T}$ is an isometry and
\[
(V_{\T}h)(\z) = \Big(\prod_{i=1}^n D_{m_i, T_i} B_{m_i}(z_i, T_i)
\Big)h \quad \quad \quad (h \in \clh, \z \in \D^n).
\]
Moreover, $V_{\T} T_i^* = M_{z_i}^* V_{\T}$, $i = 1, \ldots, n$, and
for each $\w \in \D^n$ and $\eta \in \cld_{\m,\T}$,
\[\Big((V_{\T} V_{\T}^*) (B_{\m}(\cdot, \w) \eta)\Big)(\z) =
\prod_{i=1}^n D_{m_i, T_i} B_{m_i}(z_i, T_i) B_{m_i}(w_i, T_i)^*
D_{m_i, T_i} \eta. \quad \quad (\z \in \D^n)\]
\end{Theorem}

\NI\textsf{Proof.} Clearly $V_1^*  V_1 = I_{\clh}$ and for each $j =
2, \ldots, n$, we have
\[(\bm{v}_{m_j, T_j|_{\cld_{\bm{m_{j-1},
T_{j-1}}}}})^* \bm{v}_{m_j, T_j|_{\cld_{\bm{m_{j-1}, T_{j-1}}}}} =
I_{{A^2_{\bm{m_j}}}(\cld_{\bm{m_j, T_j}})},
\]
from which we immediately deduce $V_j^* V_j =
I_{A^2_{\bm{m_{j-1}}}(\cld_{\bm{m_{j-1}}, \bm{T_{j-1}}})}$, and
finally $V_{\T}^* V_{\T} = I_{\clh}$. Now by (\ref{V_j}), we have
\[
\begin{split}
V_{\T} h & = V_n \cdots V_2 (V_1 h)= V_n \cdots V_3(V_2 D_{m_1, T_1}
B_{m_1}(z_1, T_1) h)\\ & =  V_n \cdots V_3(D_{m_2, T_2} B_{m_2}(z_2,
T_2) D_{m_1, T_1} B_{m_1}(z_1, T_1) h)\\ & = V_n \cdots V_3(D_{m_1,
T_1} D_{m_2, T_2} B_{m_1}(z_1, T_1) B_{m_2}(z_2, T_2) h),
\end{split}
\]
for all $h \in \clh$. Continuing this way we have
\[
(V_{\T} h)(\z) = \prod_{i=1}^n D_{m_i, T_i} B_{m_i}(z_i, T_i)h\quad
\quad \quad (h \in \clh, \z \in \D^n).
\]
This and a direct computation (or see \cite{AEM} or \cite{CV})
readily implies the intertwining property of $V_{\T}$ and
\[V_{\T}^*
(B_{\m}(\cdot, \w) \eta) = \prod_{i=1}^n B_{m_i}(w_i, T_i)^* D_{m_i,
T_i} \eta \quad \quad \quad (\w \in \D^n, \eta \in \cld_{\m, \T}).
\]
which in turn yields the last part of the theorem. This completes
the proof. \qed

The above theorem is a doubly commuting version and a particular
case of Theorem 3.16 in \cite{CV} by Curto and Vasilescu and
Corollary 16 in \cite{AEM} by Ambrozie, Englis and Muller (see also
\cite{BNS}). However, the present approach is based on the idea of
``simple tensor products of one variable dilation maps''. Moreover,
our construction of explicit dilation map is especially useful in
analytic model theory (see Section \ref{S-AM}).

Recall that a pair of commuting tuples $\T = (T_1, \ldots, T_n)$ on
$\clh$ and $\bm{S} = (S_1, \ldots, S_n)$ on $\clk$ is said to be
jointly unitarily equivalent, also denoted by $\T \cong \bm{S}$, if
there exists a unitary map $U : \clh \raro \clk$ such that $U T_i =
S_i U$, $i = 1, \ldots, n$. The following dilation result is an easy
consequence of Theorem \ref{dil1}.

\begin{Theorem}\label{dil-H}
Let $\T$ be a doubly commuting $B_{\m}$-contractive tuple on $\clh$.
Then there exists a joint $(M_{z_1}^*, \ldots, M_{z_n}^*)$-invariant
closed subspace $\clq \subseteq A^2_{\m}(\cld_{\DS})$ such that
\[(T_1, \ldots, T_n) \cong (P_{\clq} M_{z_1}|_{\clq}, \ldots,
P_{\clq} M_{z_n}|_{\clq}).\]
\end{Theorem}
\NI\textsf{Proof.} Let $\clq = \mbox{ran~} V_{\T}$, where $V_{\T}$
is the dilation map of $T$ as in Theorem \ref{dil1}. Then $\clq$ is
a joint $(M_{z_1}^*, \ldots, M_{z_n}^*)$-invariant subspace of
$A^2_{\m}(\clh)$ and $U_{\T} := V_{\T} : \clh \raro \clq$ is a
unitary map. Moreover, \[U_{\T} T_j^* = V_{\T} T_j^* = M_{z_j}^*
V_{\T} = M_{z_j}^* V_{\T} V_{\T}^* V_{\T} = (M_{z_j}^*|_{\clq})
V_{\T} = (M_{z_j}^*|_{\clq}) U_{\T},\]for all $j = 1, \ldots, n$.
Hence
\[U_{\T} T_j = P_{\clq} M_{z_j}|_{\clq} U_{\T} \quad \quad \quad (j = 1, \ldots,
n).\]This completes the proof of the theorem. \qed

\newsection{Analytic model}\label{S-AM}

We begin with the following lemma, the relevance of which to our
purpose will become apparent in connection with the analytic model
of doubly commuting tuples of operators.

\begin{Lemma}\label{RD}
Let $\T$ be a doubly commuting $B_{\m}$-contractive tuple on $\clh$
and $1 \leq j \leq n$. Then $A^2_{\m}(\cld_{\m, \T}) \subseteq
A^2_{\m}(\cld_{m_j, T_j})$ and $A^2_{\m}(\cld_{\m, \T})$ is a
reducing subspace for
\[
\Big(\bigotimes_{\substack{{i=1}\\i \neq j}}^n I_{A^2_{m_i}}\Big)
\otimes \bv_{m_j, T_j} \bv_{m_j, T_j}^*  \in
\clb(A^2_{\m}(\cld_{m_j, T_j})).
\]
\end{Lemma}

\NI\textsf{Proof.}  The first part follows from the fact that
$\cld_{\m, \T} \subseteq\cld_{m_j, T_j}$. For the second part, it is
enough to prove that $X_j (B_{\m}(\cdot, \w) \eta ) \in
A^2_{\m}(\cld_{\m, \T})$ where $\w \in \D^n$, $\eta \in \cld_{\m,
\T}$ and
\[ X_j :=
\Big(\bigotimes_{\substack{{i=1}\\i \neq j}}^n I_{A^2_{m_i}}\Big)
\otimes \bv_{m_j, T_j} \bv_{m_j, T_j}^*.
\]
To this end, for each $z_j, w_j \in \D$, we compute
\[
\begin{split}
D_{m_j, T_j} B_{m_j} & (z_j, T_j) B_{m_j}(w_j,
T_j)^* D_{m_j, T_j} \Big(\prod_{i=1}^n D_{m_i, T_i}\Big) \\ & =
D_{m_j, T_j} B_{m_j}(z_j, T_j) B_{m_j}(w_j, T_j)^*
\Big(\prod_{\substack{{i=1}\\i \neq j}}^n D_{m_i, T_i}\Big)
D^2_{m_j, T_j} \\ & =  D_{m_j, T_j} \Big(\prod_{\substack{{i=1}\\i
\neq j}}^n D_{m_i, T_i}\Big) B_{m_j}(z_j, T_j) B_{m_j}(w_j, T_j)^*
D^2_{m_j, T_j} \\ & = \Big(\prod_{\substack{{i=1}}}^n D_{m_i,
T_i}\Big) B_{m_j}(z_j, T_j) B_{m_j}(w_j, T_j)^* D^2_{m_j, T_j}.
\end{split}
\]
In particular, we have
\[
\Big(D_{m_j, T_j} B_{m_j} (z_j, T_j) B_{m_j}(w_j, T_j)^* D_{m_j,
T_j}\Big) \cld_{\m, \T} \subseteq \cld_{\m, \T} \quad \quad (z_j,
w_j \in \D).
\]
If $\w \in \D^n$, $\eta = \prod_{i=1}^n D_{m_i, T_i} h \in \DS$ and
$h \in \clh$, then (\ref{vv*}) gives
\[
X_j (B_{\m}(\cdot, \w) \eta ) = \Big(\prod_{\substack{{k=1}\\k \neq
j}}^n B_{m_k}(z_k, w_k)\Big) \Big(D_{m_j, T_j} B_{m_j}(z_j, T_j)
B_{m_j}(w_j, T_j)^* D_{m_j, T_j} \eta\Big) \in A^2_{\m}(\DS).
\]
This proves the desired claim, and the result follows. \qed

Therefore, $R_j \in \clb(A^2_{\m}(\DS))$, $j = 1, \ldots, n$, where
\[
R_j = \Big(\bigotimes_{\substack{{i=1}\\i \neq j}}^n
I_{A^2_{m_i}}\Big) \otimes (\bv_{m_j, T_j} \bv_{m_j,
T_j}^*)|_{A^2_{m_j}(\DS)}.
\]
By virtue of (\ref{vv*}) we have in particular for $\w \in \D^n$,
$\eta \in {\DS}$, and for $j = 1, \ldots, n$:
\begin{equation}\label{R-def}
R_j (B_{\m}(\cdot, \w) \eta) = \Big(\prod_{\substack{{i=1}\\i \neq
j}}^n B_{m_i}(z_i, w_i)\Big) \Big( D_{m_j, T_j} B_{m_j}(z_j, T_j)
B_{m_j}(w_j, T_j)^* D_{m_j, T_j} \eta \Big).
\end{equation}

\NI\textsf{Claim:} $\{R_1, \ldots, R_n\}$ is a family of commuting
orthogonal projections.

\NI\textsf{Proof of the claim:} Since $\bv_{m_j, T_j}$ is an
isometry, we deduce from the definition of $R_j$ that $R_j = R_j^* =
R_j^2$, $1 \leq j \leq n$, that is, $\{R_j\}_{j=1}^n$ is a family of
orthogonal projections. Now let $p, q \in \N$, $p \neq q$, $\w \in
\D^n$ and $\eta \in \DS$. Using (\ref{R-def}), we obtain
\[
\begin{split}
R_p R_q (B_{\m}(\cdot, \w) \eta) = & R_p
\Big(\prod_{\substack{{i=1}\\i \neq q}}^n B_{m_i}(z_i, w_i) (
D_{m_q, T_q} B_{m_q}(z_q, T_q) B_{m_q}(w_q, T_q)^* D_{m_q, T_q}
\eta) \Big) \\& = \prod_{\substack{{i=1}\\i \neq p, q}}^n
B_{m_i}(z_i, w_i) (D_{m_p, T_p} B_{m_p}(z_p, T_p) B_{m_p}(w_p,
T_p)^* D_{m_p, T_p}) \\ & \quad \quad(D_{m_q, T_q} B_{m_q}(z_q, T_q)
B_{m_q}(w_q, T_q)^* D_{m_q, T_q} \eta) \\& =
\prod_{\substack{{i=1}\\i \neq p, q}}^n B_{m_i}(z_i, w_i) (D_{m_p,
T_p} D_{m_q, T_q} B_{m_p}(z_p, T_p) B_{m_q}(z_q, T_q)
\\ & \quad \quad  B_{m_p}(w_p, T_p)^* B_{m_q}(w_q, T_q)^* D_{m_p,
T_p}D_{m_q, T_q})\eta \\& = R_q R_p (B_{\m}(\cdot, \w) \eta).
\end{split}
\]
Therefore $R_p R_q = R_q R_p$ for all $p, q = 1, \ldots, n$. The
proof of the claim is now complete.

We turn now to investigate the product $\prod_{j=1}^n R_j$. For sake
of computational simplicity, let us assume, for each $i = 1, \ldots,
n$,
\[
f_i(z, w) := B_{m_i}(z, T_i) B_{m_i}(w, T_i)^* \quad \quad \quad (z,
w \in \D).
\]
For each $\w \in \D^n$ and $\eta \in \DS$, we have
\[
\begin{split}
(\prod_{j=1}^n R_j)(B_{\m}(\cdot, \w) \eta) & =
\prod_{\substack{{j=1}\\j \neq 1}}^n R_j \Big(R_1 B_{\m}(\cdot, \w)
\eta\Big) \\ & = \prod_{\substack{{j=1}\\j \neq 1}}^n R_j \Big(
\prod_{\substack{{i=1}\\i \neq 1}}^n B_{m_i}(\cdot, w_i) D_{m_1,
T_1} f_1(\cdot, w_1) D_{m_1, T_1} \eta\Big) \\ & =
\prod_{\substack{{j=1}\\j \neq 1, 2}}^n R_j
\Big(\prod_{\substack{{i=1}\\i \neq 1, 2}}^n B_{m_i}(\cdot, w_i)
D_{m_1, T_1} D_{m_2, T_2} f_1(\cdot, w_1) f_2(\cdot, w_2) D_{m_1,
T_1} D_{m_2, T_2}\eta\Big).
\end{split}
\]
Continuing this way we have
\[
(\prod_{j=1}^n R_j)(B_{\m}(\cdot, \w) \eta)  = \prod_{i=1}^n D_{m_i,
T_i} f(\cdot, w_i) D_{m_i, T_i}  = \prod_{i=1}^n D_{m_i, T_i}
B_{m_1}(\cdot, T_i) B_{m_i}(w_i, T_i)^* D_{m_i, T_i},
\]
and hence Theorem \ref{dil1} yields $V_{\T} V_{\T}^* = \prod_{i=1}^n
R_i$. Summing up, we obtain the following:

\begin{Theorem}\label{VR}
Let $\T$ be a doubly commuting $B_{\m}$-contractive tuple on $\clh$.
Then $\{R_i\}_{i=1}^n$ is a family of commuting orthogonal
projections and  \[V_{\T} V_{\T}^* = \prod_{i=1}^n R_i.\]
\end{Theorem}

We need to introduce one more notion. For $\m \in \mathbb{N}^n$ and
for each $j = 1, \ldots, n$, set
\[
\hat{\m}_j : = (m_1, \ldots, m_{j-1},
\underbrace{1}\limits_{j-\textup{th}}, m_{j+1}, \ldots, m_n).
\]
In particular,
\[
A^2_{\hat{\m}_j} = A^2_{m_1} \otimes \cdots \otimes A^2_{m_{j-1}}
\otimes H^2 \otimes A^2_{m_{j+1}} \otimes \cdots \otimes A^2_{m_n}.
\]

Now let $\T$ be a doubly commuting $B_{\m}$-contractive tuple on
$\clh$ and $1 \leq j \leq n$. Then
\[
\begin{split}
\mbox{ran} {R}_j  & = \Big(\Big( \bigotimes_{\substack{{i=1}\\i \neq
j}}^n A^2_{m_i} \Big) \otimes \mbox{ran} (\bv_{m_j, T_j} \bv_{m_j,
T_j}^*)\Big) \bigcap \Big(\Big( \bigotimes_{\substack{{i=1}\\i \neq
j}}^n A^2_{m_i} \Big) \otimes A^2_{m_j}(\cld_{\m,\T})\Big)\\ & =
\Big( \bigotimes_{\substack{{i=1}\\i \neq j}}^n A^2_{m_i} \Big)
\otimes {\clq}_j,
\end{split}
\]
where $\clq_j = \mbox{ran} (\bv_{m_j, T_j} \bv_{m_j, T_j}^*) \cap
A^2_{m_j}(\DS)$. But $\mbox{ran} (\bv_{m_j, T_j} \bv_{m_j, T_j}^*)$
is an $M_z^*$-invariant subspace of $A^2_{m_j}(\cld_{m_j, T_j})$,
hence ${\clq}_j$ is an $M_z^*$-invariant closed subspace of
$A^2_{m_j}(\DS)$. By Theorem \ref{MT-blh} there exists an auxiliary
Hilbert space ${\cle}_j$ and a partially isometric multiplier
${\theta}_j \in \clm(A^2_1(\cle_j), A^2_{m_j}(\cld_{\m, \T}))$ such
that
\[
{R}_j = I_{A^2_{\m}(\cld_{\m, \T})} - M_{{\Theta}_j}
M_{{\Theta}_j}^*,
\]
where
\begin{equation}\label{TT}
M_{\Theta_j} = \Big(\bigotimes_{\substack{{i=1}\\i \neq j}}^n
I_{A^2_{m_i}}\Big) \otimes M_{\theta_j}.
\end{equation}
Notice that $\Theta_j \in \clm(A^2_{\hat{\m}_j}(\cle_j),
A^2_{\m}(\DS))$ is a partially isometric multiplier and
\[
\Theta_j(\z) = \theta_j(z_j) \quad \quad \quad(\z \in \D^n).
\]

This and Theorem \ref{VR} yields the following:

\begin{Theorem}\label{VTheta}
Let $\T$ be a doubly commuting $B_{\m}$-contractive tuple on $\clh$.
Then there exists Hilbert spaces $\cle_k$ and partially isometric
multipliers $\theta_k \in \clm(A^2_1(\cle_k), A^2_{m_k}(\DS))$, $k =
1, \ldots, n$, such that
\[(M_{\Theta_i} M_{\Theta_i}^*) (M_{\Theta_j} M_{\Theta_j}^*) =
(M_{\Theta_j} M_{\Theta_j}^*) (M_{\Theta_i} M_{\Theta_i}^*),\]for
all $i, j, k = 1, \ldots, n$. Moreover
\[
V_{\T} V_{\T}^* = \prod_{i=1}^n (I_{A^2_{\m}(\DS)} - M_{\Theta_i}
M_{\Theta_i}^*).
\]
\end{Theorem}

In order to formulate our functional model for $B_{\m}$-contractive
tuples, we need to recall the following result concerning commuting
orthogonal projections (cf. Lemma 1.5 in \cite{S-JOT}):

\begin{Lemma}\label{P-F} Let $\{P_i\}_{i=1}^n$ be a collection of commuting orthogonal
projections on a Hilbert space $\clh$. Then $\cll :=
\mathop{\sum}_{i=1}^n \mbox{ran} P_i$ is closed and the orthogonal
projection of $\clh$ onto $\cll$ is given by $P_{\cll} = I -
\mathop{\prod}_{i=1}^n (I - P_i)$.
\end{Lemma}

We are now ready to present the main theorem of this section.

\begin{Theorem}\label{MT}
Let $\T$ be a doubly commuting $B_{\m}$-contractive tuple on $\clh$.
Then there exists Hilbert spaces $\cle_k$ and partial isometric
multipliers $\theta_k \in \clm(A^2_1(\cle_k), A^2_{m_k}(\DS))$, $k =
1, \ldots, n$, such that
\[
\clh \cong \clq_{\Theta} := A^2_{\m}(\DS)/ \sum_{i=1}^n \Theta_i
A^2_{\hat{\m}_i}(\cle_i),
\]
and
\[
(T_1, \ldots, T_n) \cong (P_{\clq_{\Theta}}
M_{z_1}|_{\clq_{\Theta}}, \ldots, P_{\clq_{\Theta}}
M_{z_n}|_{\clq_{\Theta}}),
\]
where $\Theta_i$ is the one variable multiplier corresponding to
$\theta_i$,  $i = 1, \ldots, n$, as defined in (\ref{TT}).
\end{Theorem}

\NI\textsf{Proof.} We continue with the notation of Theorem
\ref{VTheta}. Set $P_i := M_{\Theta_i} M_{\Theta_i}^*$, $i = 1,
\ldots, n$. By virtue of Theorem \ref{VTheta} we have
\[
I_{A^2_{\m}(\DS)} - V_{\T} V_{\T}^* = I_{A^2_{\m}(\DS)} -
\prod_{i=1}^n (I_{A^2_{\m}(\DS)} - P_i).
\]
Now by Lemma \ref{P-F},
it follows that
\[
(\mbox{ran} V_{\T})^\perp = \sum_{i=1}^n \mbox{ran} M_{\Theta_i} =
\sum_{i=1}^n \Theta_i A^2_{\hat{\m}_i}(\cle_i).
\]
Therefore,
\[
\clq_{\Theta} := \mbox{ran} V_{\T} = \Big(\sum_{i=1}^n \Theta_i
A^2_{\hat{\m}_i}(\cle_i)\Big)^{\perp} \cong A^2_{\m}(\cld_{\DS})/
\sum_{i=1}^n \Theta_i A^2_{\hat{\m}_i}(\cle_i).
\]
Now using the line of argument from the proof of Theorem \ref{dil-H}
one can prove that $(T_1, \ldots, T_n) \cong (P_{\clq_{\Theta}}
M_{z_1}|_{\clq_{\Theta}}, \ldots, P_{\clq_{\Theta}}
M_{z_n}|_{\clq_{\Theta}})$. This concludes the proof. \qed

In the special case that $\m = (1, \ldots, 1)$ we recover the
functional model for doubly commuting tuples of pure contractions
\cite{BNS}. Moreover, the methods used here are different from those
used in \cite{BNS}.

\newsection{Quotient modules of $A^2_{\m}$}\label{qm}

We have a particular interest in tuples of operators $(M_{z_1},
\ldots, M_{z_n})$ compressed to joint $(M_{z_1}^*, \ldots,
M_{z_n}^*)$-invariant subspaces of reproducing kernel Hilbert spaces
over $\D^n$. Let $\clq$ be a joint $(M_{z_1}^*, \ldots,
M_{z_n}^*)$-invariant closed subspace of $A^2_{\m}$ and $C_{z_i} =
P_{\clq} M_{z_i}|_{\clq}$, $i= 1, \ldots, n$. Then $\clq$ is called
a doubly commuting quotient module of $A^2_{\m}$ if
\[
C_{z_i}^* C_{z_j} - C_{z_j} C_{z_i}^* = 0 \quad \quad (1 \leq i < j
\leq n).
\]
First, we compute the defect operator of a given doubly commuting
quotient module $\clq$ of $A^2_{\m}$:
\[
D^2_{\m, \bm{C_z}} = \prod_{i=1}^n B_{m_i}^{-1}(z_i, w_i) (C_z,
C_z^*) = P_{\clq} \Big(\prod_{i=1}^n B_{m_i}^{-1}(M_{z_i},
M_{z_i}^*)\Big)|_{\clq}.
\]
On the other hand, it is easy to see that (cf. Theorem 3.3 in
\cite{CDS})
\[
D^2_{\m, \bm{M_z}} = \prod_{i=1}^n B_{m_i}^{-1}(z_i, w_i) (M_z,
M_z^*) = \prod_{i=1}^n B_{m_i}^{-1}(M_{z_i}, M_{z_i}^*) =
P_{\mathbb{C}},
\]
where $P_{\mathbb{C}}$ is the orthogonal projection of $A^2_{\m}$
onto the one dimensional subspace of all constant functions.
Consequently, $D^2_{\m, \bm{C_z}} = P_{\clq}
P_{\mathbb{C}}|_{\clq}$, and hence
\begin{equation}\label{C}
\mbox{rank} D_{\m, \bm{C_z}} \leq1.
\end{equation}

\begin{Theorem}
Let $\clq$ be a quotient module of $A^2_{\m}$. Then the following
conditions are equivalent:

(i) $\clq$ is doubly commuting.

(ii) There exists $M_z^*$-invariant closed subspace $\clq_i$ of
$A^2_{m_i}$, $i = 1, \ldots, n$, such that \[\clq \cong \clq_1
\otimes \cdots \otimes \clq_n.\]

(iii) There exists Hilbert spaces $\cle_i$ and partially isometric
multipliers $\theta_i \in \clm(A_1(\cle_i), A^2_{m_i})$, $i = 1,
\ldots, n$, such that
\[
\clq \cong \clq_{\theta_1} \otimes \cdots \otimes
\clq_{\theta_n},
\]
where $\clq_{\theta_j} = A^2_{m_j}/ \theta_j A_1(\cle_j)$, $j = 1,
\ldots, n$.
\end{Theorem}
\NI\textsf{Proof.} Let us begin by observing that the representation
of $C_{z_i}$, $i = 1, \ldots, n$, on $\clq = \clq_1 \otimes \cdots
\otimes \clq_n$ is given by
\[
\begin{split}
C_{z_i} & = P_{\clq} (I_{A^2_{m_1}} \otimes \cdots \otimes I_{A^2_{m_{i-1}}}
\otimes M_z \otimes I_{A^2_{m_{i+1}}} \otimes \cdots \otimes
I_{\clq_n})|_{\clq} \\ & = I_{\clq_1} \otimes \cdots \otimes
I_{\clq_{i-1}} \otimes P_{\clq_i} M_z|_{\clq_i} \otimes
I_{\clq_{i+1}} \otimes \cdots \otimes I_{\clq_n}.
\end{split}
\]
This yields $(ii) \Raro (i)$ and $(iii) \Raro (i)$. The implication
$(ii) \Raro (iii)$ follows from Theorem \ref{MT-blh} and $(iii)
\Raro (ii)$ is trivial. Hence it suffices to show $(i) \Rightarrow
(iii)$. Assume $(i)$. Then by Theorem \ref{MT}, there exists Hilbert
spaces $\cle_i$ and one variable partially isometric multipliers
$\Theta_i \in \clm(A^2_{\hat{\m}_i}(\cle_i), A^2_{\m}(\cld_{\m,
\bm{C_z}}))$, $i = 1,\ldots, n$, such that
\[
(C_{z_1}, \ldots,
C_{z_n}) \mbox{~on~} \clq \cong (P_{\clq_{\Theta}}
M_{z_1}|_{\clq_{\Theta}}, \ldots, P_{\clq_{\Theta}}
M_{z_n}|_{\clq_{\Theta}}) \mbox{~on~} \clq_{\Theta},
\]
and
\[
\clq_{\Theta} = A^2_{\m}(\cld_{\m, \bm{C_z}})/\sum_{i=1}^n \Theta_i
A^2_{\hat{\m}_i}(\cle_i).
\]
Now by virtue of (\ref{C}) we have $\cld_{\m, \bm{C_z}} \cong
\{0\}$, or $\mathbb{C}$. In order to avoid trivial considerations we
assume that $\cld_{\m,\bm{C_z}} \cong \mathbb{C}$. Then
\[
\clq \cong \clq_{\Theta} = A^2_{\m}/ \sum_{i=1}^n \Theta_i
A^2_{\hat{\m}_i}(\cle_i).
\]
In particular,
\[
P_{\clq_{\Theta}} =\prod_{i=1}^n (I_{A^2_{\m}} - M_{\Theta_i} M_{\Theta_i}^*) =
\bigotimes_{i=1}^n (I_{A^2_{{m}_i}} - M_{\theta_i} M_{\theta_i}^*),
\]
which implies $\clq \cong \clq_{\Theta} = \clq_{\theta_1} \otimes
\cdots \otimes \clq_{\theta_n}$ and concludes the proof. \qed

The implication $(i) \Raro (ii)$ in previous theorem was obtained in
\cite{CDS}. For the Hardy space case $H^2(\mathbb{D}^n)$, that is,
for the case $\m = (1, \ldots, 1)$, this was observed in \cite{BNS}
and \cite{S-JOT}. Moreover, as we shall see in the next section, the
same result holds for more general reproducing kernel Hilbert spaces
over $\D^n$.

\newsection{$\frac{1}{K}$-calculus and $K$-contractivity}\label{K}

The key concept in our approach is the natural connections between
(i) operator positivity, implemented by the inverse of a positive
definite kernel function on $\D$, and a dilation map, again in terms
of the kernel function, (ii) tensor product structure of reproducing
kernel Hilbert spaces on $\D^n$, and (iii) operator positivity,
implemented by the product of $n$ positive definite kernel functions
on $\D$, of doubly commuting $n$ tuple of operators. Consequently,
our considerations can be applied even for a more general framework
(in the sense of Arazy and Englis \cite{AE}).

Let $k$ be a positive definite kernel function on $\mathbb{D}$ and
that $k(z, w)$ is holomorphic in $z$ and anti-holomorphic in $w$,
and $k(z, w) \neq 0$ for all $z, w \in \D$. Let $\clr_k \subseteq
\clo(\D, \C)$ be the corresponding reproducing kernel Hilbert space.
Moreover, let

(i) $\C[z]$ is dense in $\clr_k$,

(ii) the multiplication operator $M_z$ on $\clr_k$ is a contraction,

(iii) there exists a sequence of polynomials $\{p_k\}_{k=0}^\infty
\in \C[z, \bar{w}]$ such that
\[p_k(z, \bar{w}) \longrightarrow \frac{1}{k(z, {w})}, \quad \quad
\quad (z, w \in \D)\]and\[\sup_{k} \|p_k(M_z, M_z^*)\| < \infty.\]
We will call such a reproducing kernel Hilbert space a
\textit{standard reproducing kernel Hilbert space}, or, just SRKH
for short.

Let $\clr_k$ be a SRKH and, by virtue of condition (i) in the above
definition, let $\{\psi_k\}_{k=0}^\infty \subseteq \C[z]$ be an
orthonormal basis of $\clr_k$. For any nonnegative operator $C$ and
a bounded linear operator $T$ on a Hilbert space $\clh$, set
\[f_{k,C}(T) = I_{\clh} - \sum_{0 \leq m < k} \psi_m(T) C
\psi_m(T)^*.\]

\begin{Definition} Let $\clr_k$ be a SRKH and $T \in \clb(\clh)$. Then $T$ is said to be $k$-contractive if $\sup_{k} \|p_k(T, T^*)\| < \infty$ and
\[
C: = WOT-\lim_{k \raro \infty} p_k(T, T^*),
\]
exists and nonnegative, and
\[SOT-\lim_{k \raro \infty} f_{k, C}(T) = 0.\]
\end{Definition}

We are now ready to state the Arazy-Englis dilation result (see
Corollary 3.2 in \cite{AE}).

\begin{Theorem}
Let $\clr_k$ be a SRKH and $T \in \clb(\clh)$ be a $k$-contraction
and $\cld = \overline{\mbox{ran}}C$. Then there exists an
$M_z^*$-invariant closed subspace $\clq$ of $\clr_k \otimes \cld$
such that $T \cong P_{\clq} M_z|_{\clq}$.
\end{Theorem}

In this case, the dilation map $V_{T}$ is given by (see the equality
(1.5) in \cite{AE}): \[(V_{T} h)(z) = \sum_k \psi_k(z) \otimes
C^{\frac{1}{2}} \psi_k(T)^* h. \quad \quad \quad (h \in \clh)\]

Finally, note that the statement in Theorem \ref{MT} can be
generalized in this framework as follows (see Theorem 2.3 in
\cite{JS}): Let $\clh$ be a Hilbert space and $\cls$ be a closed
subspace of $\clr_k \otimes \clh$. Then $\cls$ is $M_z$-invariant if
and only if $\cls = \Theta H^2(\cle)$ for some Hilbert space $\cle$
and partially isometric multiplier $\Theta \in \clm(H^2(\cle),
\clr_k \otimes \clh)$.

\NI Consequently, Theorem \ref{model-1} holds for the class of
$k$-contractions.

\begin{Theorem}
Let $\clr_k$ be a SRKH and $T \in \clb(\clh)$ be a $k$-contraction
and $\cld = \overline{\mbox{ran}}C$. Then there exists a Hilbert
space $\cle$ and a partially isometric multiplier $\theta \in
\clm(H^2(\cle), \clr_k \otimes \cld)$ such that $T \cong
P_{\clq_{\theta}} M_z|_{\clq_{\theta}}$ where $\clq_{\theta} =
(\clr_k \otimes \cld) \ominus \theta H^2(\cle)$.
\end{Theorem}

Now let $\clr_{k_i}$, $i = 1, \ldots, n$, be $n$ standard
reproducing kernel Hilbert spaces over $\D$ and let
\[\clr_K := \clr_{k_1} \otimes \cdots \otimes \clr_{k_n}.\]Then
$\clr_K$ is a reproducing kernel Hilbert space (see Tomerlin
\cite{T}) and
\[
K(\z, \w) = \prod_{i=1}^n k_i(z_i, w_i) \quad  \quad \quad (\z, \w
\in \D^n).
\]

Let $\T$ be a doubly commuting tuple of operators on $\clh$ and let
$T_i$ be a $k_i$-contraction, $i = 1, \ldots, n$. Set
\[C_i = WOT-\lim_{k \raro \infty} p_{i,k}(T_i, T_i^*),\]where
$p_{i,k}(z, \bar{w}) \raro \frac{1}{k_i(z, w)}$, $i = 1, \ldots, n$.
In a similar way, as in Lemma \ref{Def1}, one can prove that $C_i
C_j = C_j C_i$ for all $i, j = 1, \ldots, n$, and
\[
C_{\T}:= \prod_{i=1}^n C_i \geq 0.
\]
By virtue of this observation, a doubly commuting tuple $\T$ is
called $K$-contractive if $T_i$ is $k_i$-contractive for all $i = 1,
\ldots, n$ (see the remark at the end of Lemma \ref{lemma1}).
Consequently, all the results and proofs in this paper hold verbatim
for this notion of a doubly commuting $K$-contractive tuples as
well.

\vspace{0.2in}

\NI\textit{Acknowledgement:} We are indebted to the referee for
numerous comments and suggestions which improved this paper
considerably. The research of the second author was supported in
part by an NBHM Research Grant NBHM/R.P.64/2014.

\end{document}